\RequirePackage[ngerman,english]{babel}
\documentclass[12pt,a4paper]{amsart}

\usepackage{amsfonts, amsmath, amssymb, amsthm}

\usepackage[ngerman,english]{babel}
\usepackage{german}
\usepackage{setspace}

\theoremstyle{plain}
\newtheorem{thm}{Theorem}[section]
\newtheorem*{thm*}{Theorem}
\newtheorem{lem}[thm]{Lemma}
\newtheorem{prop}[thm]{Proposition}
\newtheorem{cor}[thm]{Corollary}

\theoremstyle{definition}

\theoremstyle{remark}

\usepackage{algorithmic,algorithm}

\usepackage{txfonts}
\usepackage{amssymb}
\usepackage{mathbbol}
\usepackage{ae}
\usepackage{bbm}
\usepackage[mathscr]{eucal} 
\usepackage{xspace}
\usepackage{url}

\usepackage{color}
\usepackage{graphicx}
\usepackage{overpic}
\usepackage{subfigure}
\usepackage{booktabs}
\usepackage{paralist}

\usepackage[a4paper,scale=0.7, marginratio={1:1, 9:10}, ignoreall]{geometry}

\numberwithin{equation}{section}
\numberwithin{figure}{section}

\DeclareMathOperator{\Vertices}{Vert}

\renewcommand\mod{\;\mathrm{mod}\:}

\newcommand\cut{\cap}

\newcommand\union{\cup}

\newcommand\cH{{\mathcal H}}

\newcommand\cM{{\mathcal M}}

\newcommand\cT{{\mathcal T}}

\newcommand\cS{{\mathcal S}}

\newcommand\NN{{\mathbb N}}
\newcommand\RR{{\mathbb R}}

\newcommand\SetOf[2]{\left\{#1\vphantom{#2}\,\right.\left|\,\vphantom{#1}#2\right\}}
\newcommand\smallSetOf[2]{\{#1\,|\,#2\}}

\newcommand\conv{\operatorname{conv}}

\newcommand\vertt{\Vertices}

\newcommand\card[1]{\left|\,#1\,\right|}

\providecommand{\subdiv}{{\Sigma}}

\providecommand{\lift}[2]{{{\mathscr{L}}_{#1}(#2)}}
\providecommand{\tightspan}[2]{{{\mathscr{T}}_{#1}(#2)}}

\providecommand{\subdivision}[2]{{\subdiv_{#1}(#2)}}

\providecommand{\Hypersimplex}[2]{{\Delta(#1,#2)}}

\renewcommand{\phi}{\varphi}

\usepackage[T1]{fontenc}
\usepackage[utf8]{inputenc}

\graphicspath{ {pictures/} }

\begin{document}

\title{The Split Decomposition of a $k$"=Dissimilarity Map}

\author{Sven Herrmann \and Vincent Moulton}
\address{School of Computing Sciences, 
University of East Anglia, Norwich, NR4 7TJ, UK}
\email{s.herrmann@uea.ac.uk, vincent.moulton@cmp.uea.ac.uk}
\date{\today}

\begin{abstract}
A \emph{$k$"=dissimilarity map} on a finite set $X$ is a function 
$D:\binom Xk \to \RR$ assigning a real value to each subset of $X$ 
with cardinality $k$, $k \ge 2$. Such functions, also sometimes
known as $k$-\emph{way dissimilarities}, $k$-\emph{way distances},
or  $k$"=\emph{semimetrics}, are of interest in many areas of 
mathematics, computer science
and classification theory, especially 2-dissimilarity 
maps (or \emph{distances}) which
are a generalisation of metrics. In this 
paper, we show how regular subdivisions
of the $k$th hypersimplex
can be used to obtain a canonical decomposition of a 
$k$"=dissimilarity map into the sum of simpler $k$"=dissimilarity 
maps arising from bipartitions or \emph{splits} 
of~$X$. In the special case $k=2$, this
is nothing other than the well-known 
\emph{split decomposition} of a distance
due to Bandelt and Dress [Adv. Math. \textbf{92}
(1992), 47--105], a decomposition 
that is commonly to construct phylogenetic trees and 
networks. Furthermore, we characterise those 
sets of splits that may occur in the resulting 
decompositions of $k$"=dissimilarity maps.
As a corollary, we also give a new proof of a 
theorem of Pachter and Speyer 
[Appl. Math. Lett. \textbf{17} (2004), 615--621] for
recovering $k$"=dissimilarity maps from trees.
\end{abstract}
\maketitle

\section{Introduction}

Throughout this paper we assume $X=\{1,\dots,n\}$, $n \ge 1$
a natural number. For $1 < k < n$, a
\emph{$k$"=dissimilarity map} on $X$ is a 
function $D:\binom Xk \to \RR$ assigning a real value to each 
subset of $X$ with cardinality $k$ (or, alternatively 
stated, a totally symmetric function $D: X^k \to \RR$). 
Such maps are of interest in many areas of mathematics, computer
science and classification theory, especially 
$2$"=dissimilarity maps (or {\em distances}), which 
are a generalisation of metrics (cf. Deza and Laurent~\cite{MR1460488}). 
Note that $3$"=dissimilarities have been investigated, for example,
in \cite{hay-72a}, \cite{jol-cal-95}
and \cite{hei-ben-97}, and arbitrary $k$"=dissimilarities in
\cite{dez-ros-00} and \cite{war-10}, under names
such as $k$-\emph{way dissimilarities}, 
$k$-\emph{way distances} and $k$"=\emph{semimetrics}.

\begin{figure}[htb]\centering
  \input{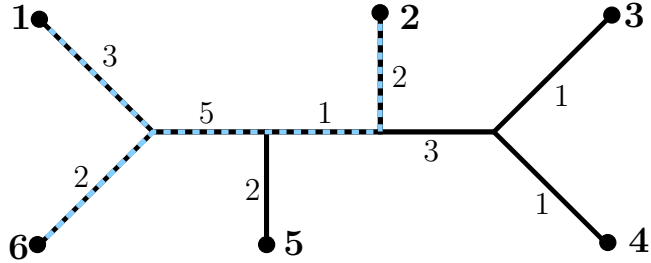}
\caption{A weighted tree, labelled by the set 
$X=\{{\bf 1},{\bf 2},\dots,{\bf 6}\}$.
A $k$"=dissimilarity map can be defined on $X$ by
assigning the length of the subtree spanned by a $k$"=subset
to that subset. For
example, if $k=3$, the subset $\{{\bf 1},{\bf 2},{\bf 6}\}$ 
would be assigned the value $13$.} 
  \label{fig:tree}
\end{figure}

Here we are interested in how to decompose $k$"=dissimilarity maps
into a sum of simpler $k$"=dissimilarity maps. Note, that 
various ways have been proposed to decompose 
distances (cf. Deza and Laurent~\cite{MR1460488}) 
although to our best knowledge not much is known 
for $k \ge 3$. More specifically, we shall introduce 
a generalisation of the {\em split decomposition} for distances
that was originally introduced by Bandelt and 
Dress~\cite{MR1153934}. The split decomposition is 
of importance in phylogenetics, where it is 
used to construct phylogenetic 
trees and networks (see e.g. Huson and Bryant~\cite{Huson06a}).
Note that $k$"=dissimilarity maps arise naturally from 
such trees (see e.g. Figure~\ref{fig:tree} and,
\cite{Levy06beyondpairwise,MR2060009}); we shall 
discuss this connection further in Section~\ref{sec:trees}.

We now explain the basic ideas underlying our results
(see Section~\ref{sec:splits} for full definitions
of the terminology that we use).
Decompositions of $k$"=dissimilarity maps
arise in the context of polyhedral 
decompositions~\cite{Triangulations} as follows.
Let $\Hypersimplex kn$ denote the
\emph{$k$th hypersimplex} $\Hypersimplex kn\subset\RR^n$,
that is, the convex hull of all $0/1$-vectors in $\RR^n$ 
having exactly $k$ ones.
Clearly, $k$"=dissimilarity maps on the set $X$ are in bijection 
with real-valued maps from the vertices of $\Hypersimplex kn$ 
since we can identify the vertices 
of $\Hypersimplex kn$ with subsets of $X$ of cardinality $k$. 
In particular, it follows that each 
$k$"=dissimilarity map~$D$ gives rise to 
a (regular) subdivision of $\Hypersimplex kn$ 
into smaller polytopes or~\emph{faces}.
We shall call a decomposition $D=D_1+D_2$ of 
$D$ \emph{coherent}, if the subdivisions 
of $\Hypersimplex kn$ corresponding 
to $D_1$ and $D_2$ have a common refinement, 
which is essentially a subdivision of
$\Hypersimplex kn$ which contains both subdivisions.

The simplest possible regular subdivision of the 
polytope $\Hypersimplex kn$ is a \emph{split subdivision} 
(or {\em split} of $\Hypersimplex kn$)~\cite{MR2502496}, 
that is, a subdivision having exactly two maximal faces.
As we shall show, using the polyhedral Split Decomposition 
Theorem~\cite[Theorem~3.10]{MR2502496}, 
it follows that a $k$"=dissimilarity 
map $D$ can always be coherently decomposed as follows.
To each bipartition or {\em split} $S=\{A,B\}$ of $X$ 
associate the {\em split $k$"=dissimilarity},
defined by
\[ \delta^k_S(K):=\begin{cases} 1,&
\text{ if }A\cap K,B\cap K\not=\emptyset,\\ 0,&
\text{else,}\end{cases} \text{for all }K\in \binom Xk\,.
\]
In addition, define the \emph{split index} 
$\alpha^D_S$ of $D$ with respect to $S$ in case $S$ is 
non"=trivial (i.e., $\card A,\card B>1$) to be the maximal 
$\lambda\in \RR_{\geq 0}$ such that 
$D=(D-\lambda \delta^k_S)+\lambda \delta^k_S$ 
is a coherent decomposition of $D$.  If $\alpha^D_S=0$ for
all splits $S$ of $X$, we call $D$ {\em split-prime}.
We prove the following:

\begin{thm}[Split Decomposition Theorem of a $k$"=Dissimilarity Map]\label{thm:splitdecomposition}~\\
  Each $k$"=dissimilarity map $D$ on $X$ has a coherent decomposition
  \begin{equation}\label{eq:split-decomposition}
  D \ = \ D_0 + \sum_{\text{$S$ split of $X$}} \alpha^D_S \delta^k_S \, ,
  \end{equation}
  where $D_0$ is split-prime. Moreover, this is unique among 
all coherent decompositions of~$D$ into a sum of split
$k$"=dissimilarities  and a split"=prime $k$"=dissimilarity map.
\end{thm}

In case $D$ is a distance (i.e., $k=2$) the 
decomposition in this theorem is
precisely the split decomposition of Bandelt and Dress~\cite{MR1153934}
mentioned above.
For such maps, it was shown in \cite[Theorem~3]{MR1153934} that the set $\cS_D$
of splits $S$ with $\alpha^D_S >0$, enjoys a special
property in that it is {\em weakly compatible}, that is, 
there do not exist (pairwise distinct) $i_0,i_1,i_2,i_3\in X$ 
and $S_1,S_2,S_3\in \cS_D$ with $S_l(i_0)=S_l(i_m)$ if 
and only if $m=l$, where 
$S(i)$ denotes the element in the split $S$ that contains~$i$.

\begin{figure}[thb]\centering
  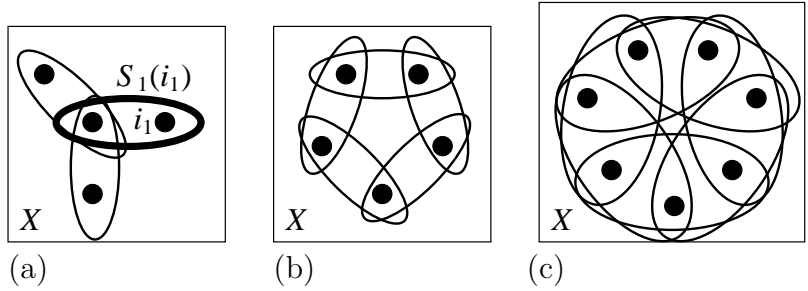
  \caption{An illustration of the forbidden situations (a)--(c)
in Theorem~\ref{thm:split-weakly-compatible}. The 
dots denote the elements $i_l\in X$ and each of the ellipses 
corresponds to one of the splits $S_l$.
For example, the dots in (a) represent the elements
$i_0,i_1,i_2,i_3$, the central dot represents the element $i_0$,
the ellipses correspond to the splits $S_1, S_2, S_3$, 
and the dots inside the bold ellipse form the set $S_1(i_1)$. The 
situations in (b) and (c) correspond to the cases 
$\nu=1$ and $\nu=7$, respectively. }
  \label{fig:forbidden}
\end{figure}

In this paper we shall show that
for a general $k$"=dissimilarity $D$, the set $\cS_D$ 
of splits with positive split index $\alpha^D_S$ can be 
characterised in a similar manner. 
In particular, calling any such set of splits 
\emph{$k$-weakly compatible}, we prove the 
following (see Figure~\ref{fig:forbidden}):

\begin{thm}\label{thm:split-weakly-compatible}
Let $\cS$ be a set of splits of  $X$. Then $\cS$ is $k$-weakly compatible if and only if none of the following conditions hold:
\begin{enumerate}
\item \label{thm:split-weakly-compatible-1} There exist (pairwise distinct) $i_0,i_1,i_2,i_3\in X$ and $S_1,S_2,S_3\in \cS$ with $S_l(i_0)=S_l(i_m) \iff m=l$
and $\big|\,{X\setminus(S_1(i_0)\cup S_2(i_0) \cup S_3(i_0))}\,\big|\geq k-2$.
\item \label{thm:split-weakly-compatible-2} For some $1\leq\nu<k$ there exist (pairwise distinct) $i_1,\dots,i_{2\nu+1}\in X$ and $S_1,\dots,S_{2\nu+1}\in \cS$ with $S_l(i_l)=S_l(i_m) \iff m\in\{l,l+1\}$ (taken modulo $2\nu+1$) and $\card{X\setminus\bigcup_{l=1}^{2\nu+1}S_l(i_l)}\geq k-\nu$.
\item \label{thm:split-weakly-compatible-3} For some $7\leq\nu<3k$ with $\nu\not\equiv 0\mod 3$ there exist (pairwise distinct) $i_1,\dots,i_{\nu}\in X$ and $S_1,\dots,S_{\nu}\in \cM$ with $S_l(i_l)=S_l(i_m) \iff m\in\{l,l+1,l+2\}$ (taken modulo $\nu$) and $\card{X\setminus\bigcup_{l=1}^{\nu} S_l(i_l)}\geq k-\lfloor\nu/3\rfloor$.
\end{enumerate}
\end{thm}

The proof of 
this characterisation will occupy a significant 
part of this paper (Section~\ref{sec:main-proof}). 
Note that it immediately follows from this theorem that 
any $k$"=weakly compatible set of splits is weakly compatible,
since the situation pictured in Figure~\ref{fig:forbidden} (a) 
is the configuration that is excluded for weakly compatible 
sets of splits in case $k=2$ (not including the 
cardinality constraint in 
Theorem~\ref{thm:split-weakly-compatible}~\eqref{thm:split-weakly-compatible-1}
which is always satisfied for $k=2$).
Also, in the special case where $D$ is a $k$"=dissimilarity map 
arising from a tree (as in \cite{hhms12}), we will further show that
Theorem~\ref{thm:splitdecomposition} can be used to recover 
the tree from~$D$ (see Theorem~\ref{thm:tree-reconstruct}). 
This gives a new proof of the main theorem 
of Pachter and Speyer in~\cite{MR2064171}.

This rest of this paper is organised as follows. 
We begin by presenting some definitions 
concerning subdivisions and splits of 
convex polytopes (Section~\ref{sec:splits}), 
as well as a short discussion on splits of 
hypersimplices (Section~\ref{sec:hyp-splits}). 
In Section~\ref{sec:split-decomposition}, we prove 
Theorem~\ref{thm:splitdecomposition}, 
while Section~\ref{sec:main-proof} is devoted to 
the rather technical proof of Theorem~\ref{thm:split-weakly-compatible}. 
This is followed by some corollaries of our main theorems 
related to $k$-weak compatibility
(Section~\ref{sec:k-weak-compat}) and 
tree reconstruction (Section~\ref{sec:trees}), respectively.
In the last section, we present some 
remarks on the connection of our results with 
tight-spans and tropical geometry as well as some 
open problems.

\noindent {\bf Acknowledgements:} The first author 
thanks the German Academic Exchange Service (DAAD) 
for its support through a fellowship within 
the Postdoc"=Programme and the UEA School of Computing Sciences
for hosting him during the writing of this paper.

\section{Subdivisions and Splits of Convex Polytopes}\label{sec:splits}

We refer the reader to Ziegler~\cite{MR1311028} and De~Loera, 
Rambau, and Santos~\cite{Triangulations} for further details 
concerning polytopes and subdivisions of polytopes, respectively. 
Let $n\geq1 $ and $P\subset\RR^{n}$ be a convex polytope. For technical reasons, we assume that $P$ has dimension $n-1$ and the origin is not an interior point of~$P$. For any hyperplane $H$ for which $P$ is entirely contained in one of the two halfspaces defined by $H$, the intersection $P\cap H$ is called a \emph{face} of~$P$.  A \emph{subdivision} of $P$ is a collection $\subdiv$ of polytopes (the \emph{faces} of $\subdiv$) such that
\begin{itemize}
\item $\bigcup_{F\in \Sigma} F=P$,
\item for all $F\in \Sigma$ all faces of $F$ are in $\Sigma$,
\item for all $F_1,F_2\in\Sigma$ the intersection $F_1\cap F_2$ is a face of $F_1$ and $F_2$,
\item for all $F\in \subdiv$ all vertices of $F$ are vertices of  $P$.
\end{itemize}

Consider a weight function $w: \Vertices P \to \RR$ assigning a weight to each vertex of~$P$. This gives rise to the {\em lifted polytope} $\lift{w}{P} \ := \conv\SetOf{(v,w(v))\in\RR^{n+1}}{v\in\Vertices{P}}$. By projecting back to the affine hull of $P$, the complex of lower faces of $\lift{w}{P}$ (with respect to the last coordinate) induces a polytopal subdivision $\subdivision{w}{P}$ of $P$. Such a subdivision of~$P$ is called a \emph{regular subdivision}. For two subdivisions $\subdiv_1,\subdiv_2$ of a polytope $P$, we can form the collection of polytopes
\begin{align}\label{eq:refinement}
\Sigma:=\smallSetOf{F_1\cap F_2}{F_1\in \subdiv_1, F_2\in \subdiv_2}\,.
\end{align}
Clearly, $\subdiv$  satisfies all but the last condition for a subdivision. If this last condition is also satisfied, the subdivision $\subdiv$ is called the \emph{common refinement} of $\subdiv_1$ and $\subdiv_2$.

A \emph{split} $S$ of $P$ is a subdivision of $P$  which
has exactly two maximal faces denoted by $S_+$ and $S_-$ (see \cite{MR2502496} for details on splits of polytopes). By our assumptions, the linear span of $S_+\cap S_-$ is a linear hyperplane $H_S$, the \emph{split hyperplane} of $S$ with respect to~$P$. Conversely, it is easily seen that a (possibly affine) hyperplane defines a split of $P$ if and only if its intersection with the (relative) interior of $P$ is nontrivial and it does not separate any edge of $P$. A set $\cT$ of splits of $P$ is called \emph{compatible} if for all $S_1,S_2\in \cT$ the intersection of $H_{S_1}\cap H_{S_2}$ with the relative interior of $P$ is empty. It is called \emph{weakly compatible} if $\cT$ has a common refinement.

\begin{lem}\label{lem:weakly-compatible-point}
Let $P$ be a polytope and $\cT$ a set of splits of $P$. Then $\cT$ is weakly compatible if and only if there does not exist a set $\cH\subset \smallSetOf{H_S}{S\in\cT}$ of splitting hyperplanes and a face $F$ of $P$ such that $F\cap \bigcap_{H\in \cH} H=\{x\}$ and $x$ is not a vertex of~$P$.
\end{lem}
\begin{proof}
Obviously, if there is a set of hyperplanes $\cH\subset \smallSetOf{H_S}{S\in\cT}$ with this property, the set $\cT$ cannot have a common refinement and hence is not compatible. Conversely, we can iteratively compute the collections \eqref{eq:refinement} for elements of $\cT$ and  it has to happen at some stage that there occurs an additional vertex $v$. At this stage take~$F$ to be the minimal face of $P$ containing $v$ and $\cH=\smallSetOf{H_S}{v\in H_S,S\in\cT}$.
\end{proof}

For a split $S$, it is easy to explicitly define a weight function $w_S$ such that $S=\subdivision{w_S}{P}$, hence all splits of~$P$ are regular subdivisions of~$P$; see \cite[Lemma~3.5]{MR2502496}. Finally, as mentioned in 
the introduction, a sum $w=w_1+w_2$ of two weight functions for $P$ is called \emph{coherent} if $\subdivision w P$ is the common refinement of $\subdivision {w_1} P$ and  $\subdivision {w_2} P$. So a sum $\sum_{S\in\cT} \lambda_S w_S$ with $\lambda_S \in \RR_{>0}$ is coherent if and only if the set $\cT$ of splits is weakly compatible.

\section{Splits of Hypersimplices}\label{sec:hyp-splits}

Let $n>k>0$. As mentioned above, the \emph{$k$th hypersimplex} $\Hypersimplex kn\subset\RR^n$ is defined as the convex hull of all $0/1$-vectors in $\RR^n$ having exactly $k$ ones, or, equivalently, $\Hypersimplex kn=[0,1]^n\cap\smallSetOf{x\in\RR^n}{\sum_{i=1}^n{x_i}=k}$. The polytope $\Hypersimplex kn$ is $(n-1)$-dimensional and has $2n$ facets defined by $x_i=1,x_i=0$ for $1\leq i\leq n$. Each face of $\Hypersimplex kn$ is isomorphic to  $\Hypersimplex{k'}{n'}$ for some $k'\leq k$, $n'<n$. This polytope first appeared in the work of Gabri\'elov, Gel$'$fand and Losik~\cite[Section~1.6]{MR0365587}.

For a split $\{A,B\}$ of $X$, and $\mu \in \NN$ the \emph{$(A,B,\mu)$-hyperplane} is defined by the equation
\begin{equation}\label{eq:hypersimplex:split}
  \mu \sum_{i\in A} x_i \ = \ (k-\mu) \sum_{i\in B} x_i \,.
\end{equation}
The splits of $\Hypersimplex kn$ can then be characterised as follows:

\begin{prop}[Lemma~5.1 and Proposition~5.2 in \cite{MR2502496}]\label{prop:hs-splits}
The splits of $\Hypersimplex kn$ are given by the $(A,B,\mu)$-hyperplanes with $k-\mu+1\le\card A\le n-\mu-1$ and $1\le\mu\le k-1$.
\end{prop}

We will be interested in the special class of splits of $\Hypersimplex kn$ defined by subsets of $X$. For $A\subsetneq X$
define the hyperplane $H_A\subset \RR^n$ by 
\begin{align}
\sum_{i\in A} x_i=1\,.
\end{align}

\begin{cor}\label{cor:H_A-splits}
For $A\subset X$ the hyperplane $H_A$ defines a split of $\Hypersimplex kn$ if and only if $2\leq \card A \leq n- k$. Otherwise, $H_A$ defines the trivial subdivision of $\Hypersimplex kn$.
\end{cor}
\begin{proof}
Since $\sum_{i=1}^n{x_i}=k$ for all $x\in \Hypersimplex kn$, the hyperplane $H_A$ defines the same split as the $(X\setminus A, A,1)$-hyperplane. Thus, by Proposition~\ref{prop:hs-splits}, $H_A$ defines a split if and only if $k\leq n-\card A\leq n-2$, which is equivalent to $2\leq \card A \leq n-k$. Obviously, if $\card A\leq 1$ or $\card A>k$, the hyperplane $H_A$ does not meet the interior of $\Hypersimplex kn$ hence defines the trivial subdivision.
\end{proof}

The split of $\Hypersimplex kn$ defined by $H_A$ for 
some $A\subset X$ will be called $S_A$. We 
now characterise when such splits of $\Hypersimplex kn$ are compatible.

\begin{lem}\label{lem:AB-compatible}
Let $A,B\subset X$. The two splits $S_A$ and $S_B$ of $\Hypersimplex kn$ are compatible if and only if either $A\subset B$, $B\subset A$, $\card{A\cup B}\geq n-k+2$, or $k=2$ and $A\cap B=\emptyset$.
\end{lem}
\begin{proof}
By \cite[Proposition~5.4]{MR2502496}, two splits of $\Hypersimplex{k}{n}$ defined by $(A,B;\mu)$- and $(C,D;\nu)$-hyperplanes are compatible if and
  only if one of the following holds:
    \begin{align*}
      \card{A\cap C} \ &\le \ k-\mu-\nu \, , &  \card{A\cap D} \ &\le \ \nu-\mu \, ,\\
      \card{B\cap  C} \ &\le \ \mu-\nu   \, , &  \text{or }\quad \card{B\cap D} \ &\le \ \mu+\nu-k \, .
    \end{align*}
That is, the two splits $S_A$ (defined by the $(X\setminus A,A,1)$-hyperplane) and $S_B$ (defined by the $(X\setminus B,B,1)$-hyperplane) are compatible if and only if 
\begin{align*}
      \card{(X\setminus A) \cap (X\setminus B)} \ &\le \ k-2 \, , &  \card{(X\setminus A)\cap B} \ &\le \ 0 \, ,\\
      \card{A\cap (X\setminus B)} \ &\le \ 0   \, , &\text{or}\mspace{60mu} \card{A\cap B} \ &\le \ 2-k \, .
    \end{align*}
The first condition can be rewritten as $\card{A\cup B}\geq n-k+2$, the second condition is equivalent to $B\subset A$, the third condition is equivalent to $A\subset B$, and the last condition can only be true if $k=2$ and $A\cap B=\emptyset$.
\end{proof}

For a weight function $w$ and a split $S_A$ of $\Hypersimplex kn$, we define the \emph{split index} $\alpha^{w}_{S_A}$ of $w$ with respect to $S_A$ as
\[ \alpha^{w}_{S_A}=\max\Big\{\lambda\in\RR_{\geq 0}\,\big|\,(w-\lambda w_{S_A})+\lambda w_{S_A} \text{ is coherent}\Big\}\,,
\]
where $w_{S_A}$ is a weight function inducing the split $S_A$ on $\Hypersimplex kn$. Note, that this is the coherency index of the weight function $w$ with respect to $w_{S_A}$ as defined in \cite[Section~2]{MR2502496}.

\section{The Split Decomposition of a $k$"=Dissimilarity Map}\label{sec:split-decomposition}

In this section, we shall prove Theorem~\ref{thm:splitdecomposition}.
We begin with some preliminaries concerning the relationship
between splits of $X$ and splits of $\Hypersimplex kn$.

As mentioned in the introduction, we can identify vertices 
of $\Hypersimplex kn$ with subsets of $X$ of cardinality $k$.  
With this identification in mind, for a $k$"=dissimilarity map~$D$, 
define the weight function 
$w_D: \vertt{\Hypersimplex kn}\to \RR;\, K\mapsto -D(K)$ on 
the vertices of~$\Hypersimplex kn$. In addition, 
for $D = \delta^k_S$, we put $w^k_S:=w^{\delta^k_S}$. This allows us to relate splits of $X$ with splits of~$\Hypersimplex kn$.

\begin{lem}\label{lem:diss-split}Let $S=\{A,B\}$ be a non"=trivial split of $X$.
\begin{enumerate}
\item \label{lem:diss-split-main}The subdivision $\subdivision{w^k_S}{\Hypersimplex kn}$ is the common refinement of the subdivisions induced on $\Hypersimplex kn$ by $H_A$ and $H_B$.
\item \begin{enumerate}
\item If $\min(\card A,\card B)\geq k$ then the subdivision $\subdivision{w^k_S}{\Hypersimplex kn}$ is the common refinement of the splits $S_A$ and $S_B$.
\item If $\card A<k\leq\card B$ then the subdivision $\subdivision{w^k_S}{\Hypersimplex kn}$ is the split $S_B$.
\item If $\max(\card A,\card B)< k$ then the subdivision $\subdivision{w^k_S}{\Hypersimplex kn}$ is trivial.
\end{enumerate}
\end{enumerate}
\end{lem}

\begin{proof}
\begin{enumerate}
\item By \cite[Lemma~3.5]{MR2502496}, a weight function for the split $S_B$ defined by the $(A,B,1)$-hyperplane is given by
\[
w_1(v)=\begin{cases} |\sum_{i=1}^n a_i v_i |, & \text{if } |\sum_{i=1}^n a_i v_i |>0,\\ 0,&\text{else},\end{cases}
\]
where $a$ is the normal vector of the $(A,B,1)$-hyperplane. Since $\sum_{i=1}^n x_i=k$ for all $x\in\Hypersimplex kn$, we have $|\sum_{i=1}^n a_ix_i |=\card {A\cap K}-(k-1)\card{B\cap K}=k(1-\card{B\cap K})$, hence (again identifying vertices of $\Hypersimplex kn$ with $k$"=subsets of $X$)
\[
w_1(K)=\begin{cases} k,& \text{if } B\cap K=\emptyset,\\ 0,&\text{else}.\end{cases}
\]
Similarly, a weight function for the split $S_A$ is given by
\[
w_2(K)=\begin{cases} k,& \text{if } A\cap K=\emptyset,\\ 0,&\text{else}.\end{cases}
\]
Obviously, $\tilde w:=\frac{w_1+w_2}k+1$ defines the same subdivision as $w_1+w_2$, and we have $\tilde w=-\delta^k_S$. 
\item Follows from \eqref{lem:diss-split-main} using Corollary~\ref{cor:H_A-splits} and Lemma~\ref{lem:AB-compatible}.
\end{enumerate}
\end{proof}

In particular, it follows from Lemma~\ref{lem:diss-split} that if $\card X\geq 2k-1$ the subdivision $\subdivision{w^k_S}{\Hypersimplex kn}$ of~$\Hypersimplex kn$ is not trivial for any split $S$, which implies in this case that the split~$S$ of~$X$ can be recovered from the subdivision $\subdivision{w^k_S}{\Hypersimplex kn}$.

Furthermore, Lemma~\ref{lem:diss-split} implies that the split index $\alpha^D_S$ of a $k$"=dissimilarity map~$D$ on~$X$ with respect to a non"=trivial split $S=\{A,B\}$ of $X$ can be written in terms of split indices for splits of the hypersimplex $\Hypersimplex kn$ as 
\[
\alpha_S^D=\min(\alpha^{w_D}_{S_B},\alpha^{w_D}_{S_A})\,.
\]
If $\alpha^D_S=0$ for all non"=trivial splits of $X$, we call $D$ \emph{free of non"=trivial splits}. This enables us to deduce our split decomposition theorem for $k$"=dissimilarities by using the polyhedral split decomposition theorem for weight functions. However, since our correspondence only works for non"=trivial splits, we have to deal with the trivial splits as a special case before we can give our proof.

\subsection{The Trivial Splits}

Each $a\in A$ defines a trivial split $S_a:=\left\{\{a\},X\setminus\{a\}\right\}$ separating $a$ from the rest of $X$. The corresponding $k$"=dissimilarity map $\delta^k_{S_a}$ on $X$ is given by
\[\delta^k_{S_a}(K):=\begin{cases} 1,&\text{if }a\in K,\\0,&\text{else.}\end{cases}
\]
Hence the extension of the weight function $w^k_{S_a}=-\delta^k_{S_a}:\Vertices \Hypersimplex kn\to \RR$ to $\RR^n$ is linear and thus induces the trivial subdivision into $\Hypersimplex kn$. In fact, $\smallSetOf{w^k_{S_a}}{a\in X}$ is a basis for the space of all functions from $\RR^n$ to $\RR$. This implies that $\alpha^{\delta^k_{S_a}}_{S}=0$ for all $a\in X$ and all non"=trivial splits $S$ of $X$, so adding or subtracting $k$"=dissimilarities corresponding to trivial splits does not interfere with split indices for non"=trivial splits.

For some $a\in X$ and a $k$"=dissimilarity map~$D$ that is free of non"=trivial splits, we define the \emph{split index} of the trivial split $S_a$ as
\[
\alpha^D_{S_a}:=\frac12\min\SetOf{\min_{b,c\in X\setminus(L\cup \{a\})}\big(D(L,a,b)+D(L,a,c)-D(L,b,c)\big)}{L\in\binom{X\setminus\{a\}}{k-2}}\,.
\]
For an arbitrary $k$"=dissimilarity map~$D$ we then set $\alpha^D_{S_a}:=\alpha^{D_0}_{S_a}$ where $D_0$ is defined as
\[D_0:=D-\sum_{\text{$S$ non"=trivial split of $X$}} \alpha^{D}_{S} \delta^k_S \,.\]

The following lemma shows that we can iteratively compute all the trivial split indices.
\begin{lem}\label{lem:trivial-unique}
Let $D$ be a $k$"=dissimilarity map on $X$, $a,a'\in X$ distinct, and $\lambda \in \RR_{\geq 0}$. Then
\[ \alpha^{D}_{S_a}=\alpha^{D+\lambda \delta^k_{S_{a'}}}_{S_a}\,.\]
\end{lem}
\begin{proof}
For all $L\in\binom{X\setminus\{a\}}{k-2}$ and $b,c\in X\setminus(L\cup \{a\})$, we see that
\begin{align*}
\delta^k_{S_{a'}}(L,a,b)+\delta^k_{S_{a'}}(L,a,c)-\delta^k_{S_{a'}}(L,b,c)&=\begin{cases} 1-1,&\text{if }a'\in L\cup\{b,c\},\\0,&\text{else,}\end{cases}\\
&=0,
\end{align*}
and hence $(D+\lambda \delta^k_{S_{a'}})(L,a,b)+(D+\lambda \delta^k_{S_{a'}})(L,a,c)-(D+\lambda \delta^k_{S_{a'}})(L,b,c)=D(L,a,b)+D(L,a,c)-D(L,b,c)$.
\end{proof}

\subsection{Proof of the Split Decomposition Theorem~\ref{thm:splitdecomposition}}

Recall that a $k$"=dissimilarity map~$D$ on $X$ is called \emph{split-prime} if for all (trivial and non"=trivial) splits~$S$ of~$X$ we have $\alpha^D_{S}=0$.

\begin{proof}
Using the Split Decomposition Theorem for polytopes~\cite[Theorem~3.10]{MR2502496}, we obtain the decomposition
\[
  w_D \ = \ w_0 + \sum_{\text{$\Sigma$ split of $\Hypersimplex kn$}} \alpha^{w_D}_{\Sigma} w_\Sigma \,,
\]
of $w_D$, where $w_\subdiv$ is a weight function defining the split $\subdiv$ of $\Hypersimplex kn$. Setting
\[
D_0:=-\left(w_0+\sum \alpha^{w_D}_{\Sigma} w_\Sigma+\sum_{A\subset X,\card A\geq2} \left(\alpha^{w_D}_{{S_A}} -\alpha^D_{\{A,X\setminus A\}}\right)w_{S_A}\right),
\]
where the first sum ranges over all splits $\Sigma$ of $\Hypersimplex kn$ that are not of the form $S_A$ for some $A\subset X$, we
can rewrite the above decomposition of $D$ as 
\begin{align*}
  D \ =  D_0 + \sum_{\text{$S$ non"=trivial split of $X$}} \alpha^{D}_{S} D^k_S \,.
\end{align*}
This decomposition is unique because of the uniqueness of the decomposition of~$w_D$.

Now for all $a\in X$ we compute the split indices  $\alpha^D_{S_a}=\alpha^{D_0}_{S_a}$ to derive the final split decomposition, which is again unique by Lemma~\ref{lem:trivial-unique}.
\end{proof}

For a $k$"=dissimilarity map $D$ on $X$, we define $\cS_D:=\smallSetOf{S\text{ split of }X}{\alpha_S^D\not=0}$, that is the set of all splits of $X$ that appear in the Split Decomposition~\eqref{eq:split-decomposition} and recall from the introduction that such a set is by definition $k$-weakly compatible.

\begin{prop}\label{prop:k-weakly-compatible}
A set $\cS$ of splits of $X$ is $k$-weakly compatible if and only if the set $\cT=\smallSetOf{S_A\text{ split of }\Hypersimplex kn}{A\in S, S\in\cS}$ of splits of $\Hypersimplex kn$ is weakly compatible.
\end{prop}
\begin{proof}
It follows from the Split Decomposition Theorem for polytopes~\cite[Theorem~3.10]{MR2502496} that a set of splits of $\Hypersimplex kn$ is weakly compatible if and only if it occurs in the split decomposition of some weight function of~$\Hypersimplex kn$. This implies that a set $\cS$ of non"=trivial splits is $k$-weakly compatible if and only if $\cT$ is a weakly compatible set of splits of $\Hypersimplex kn$. By definition, adding trivial splits does not change the $k$-weakly compatibility of a set, so the claim follows.
\end{proof}

\section{Weak compatibility of $\Hypersimplex kn$-splits}\label{sec:main-proof}

In this section, we prove a theorem from which 
Theorem~\ref{thm:split-weakly-compatible} immediately follows by Proposition~\ref{prop:k-weakly-compatible}.
For a family $\cM$ of subsets of $X$, we denote 
by $\cT(\cM):=\smallSetOf{S_A\text{ split of }\Hypersimplex kn}{A\in \cM}$ the 
corresponding set of splits of $\Hypersimplex kn$.

\begin{thm}\label{thm:hs-weakly-compatible}
Let $\cM$ be a collection of subsets of a set $X$. Then the set $\cT(\cM)$ of splits of $\Hypersimplex kn$ is weakly compatible if and only if none of the following conditions hold:
\begin{enumerate}
\item \label{thm:hs-weakly-compatible-1} There exist (pairwise distinct) $i_0,i_1,i_2,i_3\in X$ and $A_1,A_2,A_3\in \cM$ with $i_m\in A_l \iff m\in\{0,l\}$ and $\big|\,X\setminus (A_1\cup A_2 \cup A_3)\,\big|\geq k-2$.
\item \label{thm:hs-weakly-compatible-2} For some $1\leq\nu<k$ there exist (pairwise distinct) $i_1,\dots,i_{2\nu+1}\in X$ and $A_1,\dots,A_{2\nu+1}\in \cM$ with $i_m\in A_l \iff (m\in\{l,l+1\}$ (taken modulo $2\nu+1$) and $\card{X\setminus \bigcup_{i=1}^{2\nu+1} A_i}\geq k-\nu$.
\item \label{thm:hs-weakly-compatible-3} For some $7\leq\nu<3k$ with $\nu\mod 3\not=0$ there exist (pairwise distinct) $i_1,\dots,i_{\nu}\in X$ and $A_1,\dots,A_{\nu}\in \cM$ with $i_m\in A_l \iff m\in\{l,l+1,l+2\}$ (taken modulo $\nu$) and $\card{X\setminus \bigcup_{i=1}^{\nu} A_i}\geq k-\lfloor\nu/3\rfloor$.
\end{enumerate}
\end{thm}

\subsection{Sufficiency of Conditions (a)--(c)}
{\noindent (a):} Suppose \eqref{thm:hs-weakly-compatible-1} holds. Choose a subset $B$ of $X\setminus (A_1\cup A_2 \cup A_3)$ with $\card B=k-2$ and consider the face $F$ of $\Hypersimplex kn$ defined by the facets $x_i=1$ for $i\in B$ and $x_i=0$ for $i\in X\setminus (B \cup \{i_0,i_1,i_2,i_3\})$. Looking at the intersection $I:=F \cap H_{A_1}  \cap H_{A_2} \cap H_{A_3}$ we have 
\[
x_{i_0}+x_{i_1}=x_{i_0}+x_{i_2}=x_{i_0}+x_{i_3}=1 \text{ and } x_{i_0}+x_{i_1}+x_{i_2}+x_{i_3}=2 \text{ for all } x\in I\,.
\]
This yields $x_{i_k}=1-x_{i_0}$ for $k\in\{1,2,3\}$ and eventually $x_{i_k}=1/2$ for all $k\in\{0,1,2,3\}$. Hence we have $I=\{x\}$ where $x\in \RR^n$ is defined via
\[
x_i=\begin{cases}
1, &\text{if } i\in B,\\
\frac12,& \text{if } i\in\{i_0,i_1,i_2,i_3\},\\
0,&\text{else.}
\end{cases}
\]
By Lemma~\ref{lem:weakly-compatible-point}, $\cT(\cM)$ is not weakly compatible.
\\[5pt]
{\noindent (b):} Suppose \eqref{thm:hs-weakly-compatible-2} holds. Choose a subset $B$ of $X\setminus \bigcup_{i=1}^{2\nu+1} A_i$ with $\card B=k-\nu$ together with some $m\in B$ and consider the face $F$ of $\Hypersimplex kn$ defined by the facets $x_i=1$ for $i\in B\setminus\{m\}$ and $x_i=0$ for $i\in X\setminus (B \cup \{i_1,\dots,i_{2\nu+1}\})$. We consider the intersection $I:=F \cap\bigcap_{i=1}^{2\nu+1} H_{A_i}$ and get $x_{i_l}+x_{i_{l+1}}=1$ for all  $x\in I$ and $1\leq l\leq 2\nu$. So $x_{i_l}=x_{i_{l+2}}$ for all $1\leq l\leq 2\nu-1$ which implies $x_{i_1}=x_{i_{2\nu+1}}$ and, since $x_{i_{2\nu+1}}+x_{i_1}=1$, we have $x_{i_l}=1/2$ for all $1\leq k\leq 2\nu+1$. Since $\sum_{i=1}^{2\nu+1} x_i+x_m=\nu$ we also get $x_m=1/2$. Hence, we have $I=\{x\}$ where $x\in \RR^n$ is defined via
\[
x_i=\begin{cases}
1, &\text{if } i\in B\setminus\{m\},\\
\frac12,& \text{if } i\in\{i_1,\dots i_{2\nu+1},m\},\\
0,&\text{else.}
\end{cases}
\]
By Lemma~\ref{lem:weakly-compatible-point}, $\cT(\cM)$ is not weakly compatible. 
\\[5pt]
{\noindent (c):} Suppose \eqref{thm:hs-weakly-compatible-3} holds. Choose a subset $B$ of $X\setminus \bigcup_{i=1}^{\nu} A_i$ with $\card B=k-\lfloor\nu/3\rfloor$ together with some $m\in B$ and consider the face $F$ of $\Hypersimplex kn$ defined by the facets $x_i=1$ for $i\in B\setminus\{m\}$ and $x_i=0$ for $i\in X\setminus (B \cup \{i_1,\dots,i_{\nu}\})$. We consider the intersection $I:=F \cap\bigcap_{i=1}^{\nu} H_{A_i}$ and get $x_{i_l}+x_{i_{l+1}}+x_{i_{l+2}}=1$ for all  $x\in I$ and $1\leq l\leq \nu$. As in Case~\eqref{thm:hs-weakly-compatible-2} we obtain $x_{i_l}=1/3$ for all $1\leq k\leq \nu$ and, since $\sum_{i=1}^{\nu} x_i+x_m=\lfloor\nu/3\rfloor$, we get $x_m=\bar\nu/2$, where $\bar\nu=\nu\mod 3$. Hence, we have $I=\{x\}$ where $x\in \RR^n$ is defined via
\[
x_i=\begin{cases}
1, &\text{if } i\in B\setminus\{m\},\\
\frac13,& \text{if } i\in\{i_1,\dots i_{\nu},m\},\\
\frac{\bar\nu}3,& \text{if } i=m,\\
0,&\text{else.}
\end{cases}
\]
By Lemma~\ref{lem:weakly-compatible-point}, $\cT(\cM)$ is not weakly compatible.
\hfill$\square$

\subsection{Necessity of Conditions (a)--(c)}
Suppose $\cT(\cM)$ is not weakly compatible and that none of \eqref{thm:hs-weakly-compatible-1} -- \eqref{thm:hs-weakly-compatible-3} hold. Then, by Lemma~\ref{lem:weakly-compatible-point}, there exists some subset $\cM'\subset \cM$ and some face $F$ of $\Hypersimplex kn$ such that $I:=F\cap \bigcap_{A\in \cM'} H_A=\{x\}$, $x$ not a vertex of~$\Hypersimplex kn$. We assume that $\cM'$ is minimal with this property and denote by $X'\subset X$ the set of coordinates not fixed to $0$ or $1$ in $F$, that is, $0<x_i<1$ if and only if $i\in X'$. For any $i\in X'$ we denote by $\cM(i):=\smallSetOf{A\in \cM'}{i\in A}$ the set of all $A\in \cM'$ containing $i$. 

We first state some simple facts for later use:

\begin{enumerate}[(F1)]
\item \label{fact:equal} For all distinct $i,j\in X'$, we have $\cM(i)\not=\cM(j)$.
\item \label{fact:contained} For all distinct $A,B\in \cM'$, we have $A\not\subset B$.
\item \label{fact:card}For all $A\in\cM'$, we have $\card{A\cut X'}\geq 2$.
\item \label{fact:two}For all $A\in \cM'$, there exists some $i\in A$ with $\card{\cM(i)}\geq 2$.
\end{enumerate}
\begin{proof}
\begin{enumerate}[(F1)]
\item Suppose there exist distinct $i,j\in X'$, with $\cM(i)=\cM(j)$. Then choose some $0<\epsilon<\min(x_i,1-x_j)$ and consider $x'\in\RR^n$ defined by
\[
x'_l=\begin{cases}
x_l-\epsilon,&\text{if } l=i,\\
x_l+\epsilon,&\text{if } l=j,\\
x_l,&\text{else.}
\end{cases}
\]
So $x\not=x'$ and $x'\in I$, a contradiction.
\item Follows from the minimality of $\cM'$.
\item Suppose $\card{A\cut X'}=\{j\}$ for some $A\in\cM'$ and $j\in X'$. Then $0<x_j<1$ but $x_i\in\{0,1\}$ for all $i\in A\setminus\{j\}$ which obviously contradicts $\sum_{i\in A}x_i=1$.
\item Let $A\in \cM'$. By F\ref{fact:card} there exist distinct $i,=j\in A$ and by F\ref{fact:equal} $\cM(i)\not=\cM(j)$. However, $A\in\cM(i)\cap \cM(j)$ so either $\cM(i)$ or $\cM(j)$ has to contain another $B\in \cM'$.
\end{enumerate}
\end{proof}

As the next step, we will show that none of the following 
conditions may be satisfied:
\begin{enumerate}[(i)]
\item \label{cond:three}There exists (pairwise distinct) $i_0,i_1,i_2,i_3\in X'$ and $A_1,A_2,A_3\in \cM'$  with $i_m\in A_l \iff m\in\{0,l\}$.
\item \label{cond:circle-odd} For some $\nu\in \NN$, there exist (pairwise distinct) $i_1$, $\dots$, $i_{2\nu+1}\in X'$ and $A_1$, $\dots$, $A_{2\nu+1}\in \cM'$ with $i_m\in A_l \iff m\in\{l,l+1\}$ (taken modulo $2\nu+1$).  
\item \label{cond:sequence-even}For some $\nu\in \NN$, there exist (pairwise distinct) $i_0,i_1,\dots,i_{2\nu+1}\in X'$ and $A_1$, $\dots$, $A_{2\nu+1}\in \cM'$ with $\cM(i_0)=\{A_1\},\cM(i_{2\nu+1})=\{A_{2\nu+1}\}$ and $\cM(i_l)=\{A_l,A_{l+1}\}$ for $1\leq l \leq 2\nu$.
\item \label{cond:circle-even}For some $\nu\in \NN$, there exist (pairwise distinct) $i_1,\dots,i_{2\nu}\in X'$ and $A_1,\dots,A_{2\nu}\in \cM'$ with $\cM(i_l)=\{A_l,A_{l+1}\}$ (taken modulo $2\nu$).
\item \label{cond:card-three} There exists some $i\in X'$ with $\card{\cM(i)}=3$.
\item \label{cond:three-four} For some $A\in\cM'$, there exist distinct $i,j\in A$ such that $\card{\cM(i)},\card{\cM(j)}\geq 4$.
\end{enumerate}
\begin{proof}
\noindent (i): Suppose this were true. Then we have $\sum_{i\in A_l\setminus\{i_0\}} x_i= 1-x_{i_0}$ for $l\in\{1,2,3\}$, hence $\sum_{i\in A_1\cup A_2 \cup A_3}x_i\leq x_{i_0}+\sum_{l=1}^{3}\sum_{i\in A_l\setminus\{i_0\}}x_i\leq 3-2x_{i_0}<3$. Since $\sum_{i\in X} x_i=k$, this implies $\sum_{i\in X\setminus (A_1\cup A_2 \cup A_3)} x_i> k-3$ and, because $x_i\in\{0,1\}$ for all $i\in   X\setminus (A_1\cup A_2 \cup A_3)$, we get $\card{X\setminus (A_1\cup A_2 \cup A_3)}\geq k-2$. So we are in situation \eqref{thm:hs-weakly-compatible-1} of the theorem, a contradiction.
\\[5pt]
\noindent (ii): For the purpose of this proof, a collection of $i_l$ and $A_l$ satisfying this condition will be called a \emph{cycle}. We set $T=\bigcup_{i=1}^{2\nu+1}A_i$, $T_1:=\smallSetOf{i_l}{1\leq l\leq 2\nu+1}$, $T_2:=T\setminus T_1$, $t:=\card T$, $t_1:=\card {T_1}$, and $t_2:=\card{T_2}$.  Cycles are partially ordered by the lexicographic ordering of the pair $(\nu,t)$. We assume without loss of generality that our cycle is minimal in the set of all cycles occurring in $\cM'$.

As base case we consider $\nu=1$ and $t\leq 5$. Each decreasing chain of cycles will eventually reach this case since $\nu\geq1$ and $t\geq2\nu+1$. Then (after a possible exchange of $A_3$ with $A_1$ or $A_2$) we can assume that $T\subset A_1\union A_2$, hence $\sum_{i\in T}x_i<2$. This implies that $\sum_{i\in X\setminus T}x_i>k-2$ and hence $n-t\geq k-1$ since $x_i\in\{0,1\}$ for all $i\in X\setminus T$. So we are in situation~\eqref{thm:hs-weakly-compatible-2} of the theorem, a contradiction.

We say that a set $A\in \cM'$ is of \emph{a-type} (with respect to some cycle $Z$) if for some $1\leq l\leq 2\nu+1$ we have $i_l\in A$, $A\subset A_l\cup A_{l+1}$, and $\card{A\cap T_2}\geq2$. The set is of \emph{b-type} (with respect to some cycle $Z$) if there exists some $i\in A\cap T_2$ and some $j\in A\cap (X'\setminus T)$. We will show that for the cycle $Z$ each set $A$ (distinct from all $A_l$) with $A\cap T\not=\emptyset$ is either of a-type or of b-type with respect to $Z$.

First consider some set $A$ (distinct from all $A_l$) with $i_l\in A$ for some $1\leq l\leq 2\nu+1$ and some $j\in A\setminus\{i_l\}$. Then $j\in T$ because otherwise $i_l,i_{l-1},i_{l+1},j$ and $A_l,A_{l+1},A$ would satisfy Condition~\eqref{cond:three} for some $j\in A\setminus T$. Furthermore, if there exists some $m\not\in\{l,l+1\}$ with $j\in A_m$, then we could form a smaller cycle. We get $j\in A_l\cup A_{l+1}$ and (using F\ref{fact:contained}) $\card{A\cap T_2}\geq2$, so $A$ is of a-type.

Now fix a minimal cycle $Z$ and consider an arbitrary set $B$ (distinct from all~$A_l$) with $B\cap T\not=\emptyset$. Suppose that $B\subset T_2$. This implies that there either exists a smaller cycle, or we have the situation that there exists some $1\leq l\leq 2\nu+1$ such that $B\subset A_{l+1}\cup A_{l-1}$ and $B\cap A_{l+1},B\cap A_{l-1}\not=\emptyset$. By the minimality of our cycle this implies $A_l\subset T_1$. However, this implies $B\cup A_l\subsetneq A_{l+1}\cup A_{l-1}$, a contradiction to $\sum_{i\in A}x_i=1$ for all $B\in \cM'$ and $x_i>0$ for all $i\in X'$. So $B$ either contains some element of $T_1$ implying $B$ is of a-type or some element from $X\setminus T$ implying $B$ is of b-type.

Now each $i\in T_1$ cannot be contained in some set of b-type by definition and can be contained in at most one set of a-type by F\ref{fact:contained}. Furthermore, each $i\in T_2$ can be contained in at most two sets of a-type or in at most one set of b-type but not both. To see this assume that $i\in A_l$ is contained in two sets $A,B$ either $A$ of a-type and $B$ of b-type or both of b-type. Then there exist $i_1\in A\setminus (B\cup A_l)$, $i_2\in B\setminus (A\cup A_l)$, and  $i_3\in A_l\setminus (B\cup A)$ such that $A,B,A_l$ and $i,i_1,i_2,i_3$ satisfy Condition~\eqref{cond:three}. For the same reason, each $i\in X'\setminus T$ can be in at most two sets of b-type.

We denote the number of sets of a-type (b-type) with respect to $Z$ by~$a$ (by~$b$). In order to uniquely define all $t$~coordinates of $x_i$ with $i\in T$, it is necessary to have at least $t$~equations involving some $x_i$ with $i\in T$, that is, $t$~sets in $\cM'$ which contain elements of $T$. By our considerations above, all such sets have to be either of a-type or of b-type or be equal to some~$A_l$ for $1\leq l\leq 2\nu+1$. Hence we get $a+b+2\nu+1\geq t$, or, equivalently (since $t=t_1+t_2=t_2+2\nu+2$),
\begin{align}\label{eq:equations} a+b\geq t_2\,.\end{align}
Furthermore, by the fact that some $j\in T_2$ can only be in one set of b-type and this holds only if it is not in some set of a-type, we have $b\leq t_2-a'$, where $a'$ is the number of elements of $T_2$ contained in some set of a-type. Together with Inequality~\eqref{eq:equations} we obtain
\begin{align}\label{eq:a-b} t_2-a\leq b\leq t_2-a'\,;\end{align}
in particular $a'\leq a$. However, since each set of a-type contains at least two elements of~$T_2$ and each element of~$T_2$ is contained in at most two sets of a-type, which implies $a'\geq a$, we have $a'=a$ and each element of~$T_2$ is contained in either one set of b-type or each element of ~$T_2$ is contained in exactly two sets of a-type. In view of the definition of the sets of a-type the former implies that there are no sets of a-type at all and the latter implies that the sets of a-type with respect to $Z$ form themselves a cycle $Z'$ together with the elements $j_1,\dots,j_{2\nu+1}\in T_2$ contained in sets of a-type with respect to~$Z$.

We first consider the latter case. Suppose without loss of generality that $j_l\in A_l$ and call the set of a-type containing $j_l$ and $j_{l+1}$ $B_l$. Then the sets $A_l$ are sets of a-type with respect to $Z'$. Hence $\bigcup_{l=1}^{2\nu+1}(A_l\cup B_l)=\smallSetOf{i_l,j_l}{1\leq l\leq 2\nu+1}$. If now $2\nu+1$ is not divisible by $3$, then we are in the situation \eqref{thm:hs-weakly-compatible-3} of the theorem, a contradiction, since $\nu\geq 6$ by our base case and $\nu<3k$ obviously holds. If $2\nu+1$ is divisible by $3$, then choose some $0<\epsilon<\min\smallSetOf{x_{i_l},x_{j_l}}{1\leq l\leq 2\nu+1}$ and consider $x'\in\RR^n$ defined by
\[
x'_l=\begin{cases}
x_l+\epsilon,&\text{if } l=i_m \text{ and }m\equiv 1 \mod 3\text{ or }l=j_m \text{ and }m\equiv 2 \mod 3,\\
x_l-\epsilon,&\text{if } l=j_m \text{ and }m\equiv 1 \mod 3\text{ or }l=i_m \text{ and }m\equiv 3 \mod 3,\\
x_l,&\text{else}.
\end{cases}
\]
Then $x\not=x'$ and $x\in I$, a contradiction.

The case remaining is $a=0$. Then Inequality~\eqref{eq:a-b} implies $b=t_2$. So $\sum_{i\in T}x_i=2\nu+1-\sum_{i\in T_1} x_i$, since each element of $T_2$ is in exactly one of the $2\nu+1$ sets $A_l$ and each element of $T_1$ in exactly two. This is equivalent to $\sum_{i\in T_1} x_i=\nu-\frac12(\sum_{i\in T_2} x_i-1)$. Define $T_3$ to be the set of all elements of $X'$ that are one set of b-type but not in $T$. There cannot be any elements of $X'$ that are in more than two sets of b-type but not in $T$ because this would satisfy Condition~\eqref{cond:three}. For some $t\in T_3$ which is in exactly one set of b-type, we get $x_t\leq 1-x_j\leq 1-1/2x_j$ for some $j\in T_2$, and for some $t\in T_3$ which is in exactly two sets of b-type, we get $x_t\leq 1-\max(x_j,x_l)\leq 1-1/2(x_j+x_l)$ for some $j,l\in T_2$. Since each $j\in T_2$ is contained in exactly one set of b-type, each $j\in T_2$ occurs exactly once, hence we get $\sum_{i\in T_3}x_i\leq \card{T_3}-\frac12 \sum_{i\in T_2}x_i$.
So
\begin{align*}
\sum_{i\in T\cup T_3}x_i&\leq \nu-\frac12\left(\sum_{i\in T_2} x_i-1\right)+\sum_{i\in T_2}x_i+\card{T_3}-\frac12\sum_{i\in T_2}x_i\\
&=\nu+\card{T_3}+\frac12\,,
\end{align*}
and $\card{X\setminus(T\cup T_3)}\geq k-\nu-\card{T_3}-1/2)$. Hence $\card{X\setminus T}\geq k-\nu$ (as it has to be an integer). So we are in the situation \eqref{thm:hs-weakly-compatible-2} of the theorem, a contradiction
\\[5pt]
\noindent (iii),(iv): Choose some $0<\epsilon<\min\smallSetOf{x_{i_l}}{l\text{ odd}}\cup\smallSetOf{1-x_{i_l}}{l\text{ even}}$ and define the point $x'\in \RR^n$ by\
\[
x'_l=\begin{cases}
x_l-\epsilon,&\text{if } l=i_j \text{ for some odd }j,\\
x_l+\epsilon,&\text{if } l=i_j \text{ for some even }j,\\
x_l,&\text{else}.
\end{cases}
\]
Obviously, $x\not=x'$ and it is easily checked that $x'\in I$, a contradiction.
\\[5pt]
\noindent (v): Suppose there exists some $i\in X'$ with $\cM(i)\geq 3$.  Since Condition~\eqref{cond:three} cannot hold, there has to exist some $B\in \cM(i)$  such that, for each $j\in B$, there exists some $B\not=C$ with $j\in C\in\cM(i)$. By F\ref{fact:contained}, there exist distinct $j_1,j_2\in B$ and $C_1,C_2\in\cM(i)$ with $j_1\in C_1$, $j_2\in C_2$ and $l_1\in C_1\setminus B$, $l_2\in C_2\setminus B$. Furthermore, we have $C_1\cap C_2=\emptyset$ because otherwise $B,C_1,C_2$ and $j_1,j_2,j_3$ for some $j_3\in C_1\cap C_2$ would satisfy Condition~\eqref{cond:circle-odd}. So for each $i\in X'$ with $\card{\cM(i)}\geq 3$ we have the situation depicted in the left of Figure~\ref{fig:three}.

\begin{figure}[htb]\centering
  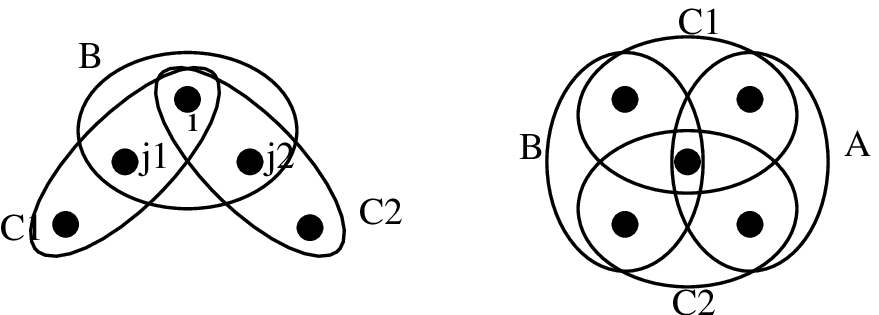
  \caption{Situations for $i\in X'$ with $\cM(i)=3$ and $\cM(i)\geq4$, respectively.}
  \label{fig:three}
\end{figure}

If there now exists some other point $i'\in C_1$ with $\card{\cM(i')}\geq 3$, then we have to be in the same situation for this point again if $i'\not\in B$. In particular this implies also that $\card{\cM(j)}\geq 3$ for some $j\in A$, so we can assume that $i'\in B$. We now repeat this process until we either get an element that we had before -- implying that Condition~\eqref{cond:circle-odd} holds -- or we arrive at some set $A$ that has exactly one $i\in A$ with $\card{\cM(i')}\geq 3$. 

Repeating the same process for $C_2$ instead of $C_1$, we finally arrive at the following situation: For some $\nu\in \NN$ there exist $i_1,\dots,i_\nu$ and $A_1,\dots,A_\nu$ such that $\cM(i_1)=\{A_1,A_2\}$, $\cM(i_l)=\{A_{l-1},A_l,A_{l+1}\}$ for $1<l<\nu$, $\cM(i_\nu)=\{A_{\nu-1},A_\nu\}$, and $A_l=\{i_{l-1},i_l,i_{l+1}\}$ for $1<l<\nu$.

We now consider two cases: First suppose $\nu\equiv 2 \mod 3$. Then choose $0<\epsilon<\min\smallSetOf{x_{i_l},1-x_{i_l}}{1\leq l\leq \nu}$ and consider $x'\in\RR^n$ defined by
\[
x'_l=\begin{cases}
x_l+\epsilon,&\text{if } l=i_m \text{ and }m\equiv 1 \mod 3,\\
x_l-\epsilon,&\text{if } l=j_m,\: l=i_m \text{ and }m\equiv 1 \mod 3,\\
x_l,&\text{else}.
\end{cases}
\]
Then $x\not=x'$ and $x\in I$, a contradiction. So suppose $\nu\not\equiv 2 \mod 3$. Then it is easily seen that the values of $x_{i_l}$ for $1\leq l\leq \nu$ are determined by the values $\sum_{i\in A_1\setminus\{i_1,i_2\}}x_i$ and $\sum_{i\in A_\nu\setminus\{i_{\nu-1},i_\nu\}}x_i$. This implies that $\cM'':=\cM'\setminus\smallSetOf{A_l}{1\leq l\leq \nu}$ with \[F':=F\cap\SetOf{x\in\RR^n}{x_{i_{l}}=\begin{cases}1,&\text{if }l\equiv1\mod 3\\0,\text{else,}\end{cases}\text{ for all }1\leq l\leq \nu}\] if $\nu\equiv 0 \mod 3$ and \[F':=F\cap\SetOf{x\in\RR^n}{x_{i_{l}}=\begin{cases}1,&\text{if }l\equiv1\mod 3\\0,\text{else,}\end{cases}\text{ for all }1\leq l< \nu}\]if $\nu\equiv 2 \mod 3$ would also have been a valid choice at the beginning, but $\cM''\subsetneq \cM'$ contradicts the minimality of $\cM'$.
\\[5pt]
\noindent (vi): Suppose there exists some $i\in X'$ with $\cM(i)\geq 4$. As in the proof of~\eqref{cond:card-three}, we have to be in the situation depicted in the left of Figure~\ref{fig:three} and there exists some $A\in\cM(i)\setminus\{B,C_1,C_2\}$. Since Condition~\eqref{cond:three} cannot hold, every $j\in A$ has to be in some $C'\in\cM(i)$ and, again by F\ref{fact:contained}, there exist distinct $j'_1,j'_2\in A$ and $C'_1,C'_2\in\cM(i)$ with $j'_1\in C_1$, $j'_2\in C'_2$ and $l'_1\in C'_1\setminus A$, $l'_2\in C'_2\setminus A$. Since Condition~\eqref{cond:three} cannot hold, we get $C'_1,C'_2\in\{A,C_1,C_2\}$. However, if, for example, $C'_1=C_1$ and $C'_2=A$, then $i,j'_1,j'_2$ and $A,B,C_1$ would satisfy Condition~\eqref{cond:circle-odd}. Hence we have  $C'_1=C_1$ and $C'_2=C_2$, or vice-versa. So we are in the situation depicted in the right of Figure~\ref{fig:three}. To obtain in addition some $j$ with $\cM(j)\leq 3$, there has to exist some $D\in\cM'$ with $D\cap U\not=\emptyset$, where $U:=A\cup B\cup C_1 \cup C_2$. Because Condition~\eqref{cond:three} cannot hold, we get $\card{D\cap U}\geq 2$ and so F\ref{fact:contained} implies that either Condition~\eqref{cond:three} or Condition~\eqref{cond:circle-odd} has to be satisfied, a contradiction.
\end{proof}

We will now show that under our assumptions at the beginning of the proof one of the Conditions \eqref{cond:three} to \eqref{cond:three-four} has to be satisfied, which leads to a contradiction.

For each $A\in \cM'$ we define $\tilde A:=\smallSetOf{i\in A\cap X'}{\cM(i)\leq 2}$. We have $\card {\tilde A}\geq 2$ for all $A\in\cM'$, because otherwise we would have a situation satisfying one of Conditions~\eqref{cond:card-three} or~\eqref{cond:three-four}. Given some pair $(A,\delta)\in\cM'\times X'$ with $\delta\in \tilde A$, we now give a way to construct a finite sequence $F(A,\delta)=(A_j,\alpha_j)_{1\leq j\leq L(A,\delta)}\subset \cM'\times X'$:
\begin{enumerate}[I]
\item $(A_1,\alpha_1):=(A,\delta)$.
\item \label{case:circle} If there exists some $\gamma\in \tilde A_j$ such that $A_l\in\cM(\gamma)$ for some $l<k$, then $L(A,\delta)=j$ and $(A_j,\alpha_j)$ is the last element of the sequence;
\item else, if there exists some $\gamma\in\tilde A_j$ such that $\cM(\gamma)=\{A_j,C\}$ for some $C\not=A_j$, then we set $A_{j+1}:=C$ and $\alpha_{j+1}:=\gamma$;
\item \label{case:chain} else, there exist a (unique) $\gamma\in\tilde A_j$ with $\cM(\gamma)=\{A_j\}$; then $L(A,\delta)=j$ and $(A_j,\alpha_j)$ is the last element of the sequence.
\end{enumerate}
The existence of the $\gamma\in \tilde A_j$ in Case~\ref{case:chain} follows from the fact that $\card{\tilde A_j} \geq 2$ and its uniqueness from F\ref{fact:equal}. Obviously, $F(A,\delta)$ ends in either Case~\ref{case:circle} or in Case~\ref{case:chain}. Suppose there exist some pair $(A,\delta)$ ending up in Case~\ref{case:circle}. Then $\alpha_1,\dots,\alpha_{L(A,\delta)}$ and $A_1,\dots,A_{L(A,\delta)}$ obviously satisfy Condition~\eqref{cond:circle-odd} if $L(A,\delta)$ is odd and Condition~\eqref{cond:circle-even} if $\mu_A$ is even -- a contradiction. Hence for each starting pair $(A,\delta)\in\cM'\times X'$ with $\delta\in \tilde A$ we end up in Case~\ref{case:chain}. The unique element $\gamma$ occurring there will be denoted $f(A,\delta)$.

Now choose some $B\in \cM'$. By F\ref{fact:two} and $\card{\tilde B}\geq 2$ there exists some $\delta\in B$ with $\card{\cM(\delta)}= 2$, say $\cM(\delta)=\{B,C\}$ for some $C\not=B$. We now construct the sequences $F(B,\delta)=(B_j,\alpha_j)_{1\leq j\leq L(B,\delta)}$ and $F(C,\delta)=(C_j,\gamma_j)_{1\leq j\leq L(C,\delta)}$. Define
\begin{align*}
i_0:=f(B,\delta),\quad i_1&:=\beta_{L(B,\delta)}&A_1&:=B_{L(B,\delta)}\\
&\dots&&\dots\\
i_{L(B,\delta))}&:=\beta_1=\delta=\gamma_1,&A_{L(B,\delta)}&:=B_1, A_{L(B,\delta)}:=B_1,\\
&\dots&&\dots\\
i_{L(B,\delta)+L(C,\delta)-1}:=\gamma_{L(C,\delta)}, i_{L(B,\delta)+L(C,\delta)}&:=f(B,\delta)&A_{L(B,\delta)+L(C,\delta)}&:=C_{L(C,\delta)}\,.
\end{align*}
Now if $e:=L(B,\delta)+L(C,\delta)$ is odd, then these $i_0,\dots,i_e$ and $A_1,\dots,A_e$ satisfy Condition~\eqref{cond:sequence-even}. So $e$ must be even.

Suppose there exists some $1<j<e$ and some $\alpha\in A_j$ with $\alpha\not=i_{j-1},i_j$. Then we distinguish two cases: First, assume that $\cM(\alpha)=\{A_j\}$. Then either $j$ is odd and $i_0,\dots,i_{j-1},\alpha$ and $A_1,\dots,A_j$ satisfy Condition~\eqref{cond:sequence-even}, or $j$ is even, hence $e-j+1$ is odd and $\alpha,i_{j}\dots,i_{e}$ and $A_{j},\dots,A_e$ satisfy Condition~\eqref{cond:sequence-even}. So assume that $D\in\cM(\alpha)$ for some $D\not=A_j$. Now we construct the sequence $F(D,\alpha)=(D_j,\delta_j)_{1\leq j\leq L(D,\alpha)}$. Then either $j+L(D,\alpha)$ is odd and $i_0,\dots,i_{j-1},\alpha=\delta_1,\dots,\delta_{L(D,\alpha)},f(D,\alpha)$ and $A_1,\dots,A_j$, $D_1,\dots,D_{L(D,\alpha)}$ satisfy Condition~\eqref{cond:sequence-even} or $j+L(D,\alpha)$ is even, hence $e-j+L(D,\alpha)+1$ is odd and, similarly, $i_e\dots,i_{j},\alpha=\delta_1,\dots,\delta_{L(D,\alpha)},f(D,\alpha)$ and $A_e,\dots,A_j,D_1,\dots,D_{L(D,\alpha)}$ satisfy Condition~\eqref{cond:sequence-even}.

This shows that for each $\alpha\in A_j$ with $1< j<e$ we have $\cM(\alpha)=\{A_{j-1},A_j\}$ or $\cM(\alpha)=\{A_{j},A_{j+1}\}$. By~F\ref{fact:equal}, this implies $\alpha=i_j$ or $\alpha=i_{j-1}$, respectively. Furthermore, it follows from this fact and the construction of $F(B,\delta)$ and $F(C,\delta)$ that $\alpha\in A_1\setminus A_2$ implies $\alpha=i_0$ and $\alpha\in A_{e}\setminus A_{e-1}$ implies $\alpha=i_e$. Thus, each $A_j$, $1\leq j\leq e$,  has exactly two elements. Hence $x$ has to satisfy the equations
\[
x_{i_l}+x_{i_{l-1}}=1,\quad\text{for all }1\leq l \leq e\,.
\]
This implies that $x_{i_l}=x_{i_0}$ if $l$ is odd and $x_{i_l}=1-x_{i_0}$ if $l$ is even. In particular, $\sum_{l=0}^{e}x_{i_l}=x_0+\sum_{l=1}^{e/2}(x_{i_l}+x_{i_{l-1}})=e/2+x_0$ is not an integer, hence there exists some $\gamma\in X'\setminus\{i_0,\dots,i_e\}$. We distinguish two cases: If $\cM(\gamma)=\emptyset$, then choose some $0<\epsilon<\min\{x_{i_0},1-x_{i_0},x_\gamma\}$ and define the point $x'\in \RR^n$ via
\[
x'_i=\begin{cases}
x_i-\epsilon,&\text{if } i=i_l \text{ for some even }l,\\
x_i+\epsilon,&\text{if } i=\gamma, \text{ or }i=i_l \text{ for some odd }l\\
x_i,&\text{else}.
\end{cases}
\]
It is easily checked that $x'\in I$, a contradiction.

In the case $\cM(\gamma)\not=\emptyset$ there exist some $B^\star\in \cM'$ with $\gamma\in B^\star$. We can now argue as before: By F\ref{fact:two} and $\card{\tilde{B^\star}}\geq 2$ there exists some $\delta^\star\in B^\star$ with $\card{\cM(\delta^\star)}=2$, say $\cM(\delta^\star)=\{B^\star, C^\star\}$ for some $C^\star\not=B^\star$. This leads us to $i^\star_0,\dots i_{e^\star}^\star$ and $A^\star_1,\dots A^\star_{e^\star}$ having the same properties as $i_0,\dots,i_e$ and $A_1,\dots,A_e$. Choose some $0<\epsilon<\min\{x_{i_0},1-x_{i_0},x_{i^\star_0},1-x_{i^\star_0}\}$ and define the point $x'\in \RR^n$ via
\[
x'_i=\begin{cases}
x_i-\epsilon,&\text{if } i=i_l \text{ for some even }l \text{ or } i=i^\star_l \text{ for some odd }l,\\
x_i+\epsilon,&\text{if } i=i_l \text{ for some odd }l \text{ or } i= i^\star_l \text{ for some even }l,\\
x_i,&\text{else}.
\end{cases}
\]
It is easily checked that $x'\in I$, our final contradiction.\hfill$\square$

\section{Compatibility and $k$"=Weak Compatibility of splits of $X$}\label{sec:k-weak-compat}

In this section, we present some corollaries of Theorem~\ref{thm:split-weakly-compatible}. Recall that two splits $\{A,B\}$ and $\{C,D\}$ are called \emph{compatible} if one of the four intersections $A\cap C$, $A\cap D$, $B\cap C$, or $B\cap D$ is empty; a set $\cS$ of splits is called \emph{compatible} if each pair of elements of $\cS$ is compatible (see e.g.,~\cite{MR2060009}).

We first consider the case $k=2$. In this case, for a split $\{A,B\}$ of $X$, the splits $S_A$ and $S_B$ of $\Hypersimplex 2n$ are clearly equal. 

\begin{cor}[Corollary~6.3 and Proposition~6.4 in~\cite{MR2502496}]\label{cor:k-2-compatible}
Let $\cS$ be a set of splits of $X$. 
\begin{enumerate}
\item $\cS$ is compatible if and only if $\cT:=\smallSetOf{S_A\text{ split of }\Hypersimplex 2n}{A\in S, S\in\cS}$ is a compatible set of splits of $\Hypersimplex 2n$
\item $\cS$ is weakly compatible if and only it is $2$-weakly compatible.
\end{enumerate}
\end{cor}
\begin{proof}
\begin{enumerate}
\item Follows from Lemma~\ref{lem:AB-compatible}.
\item Condition~\eqref{thm:split-weakly-compatible-1} of Theorem~\ref{thm:split-weakly-compatible} reduces exactly to the usual definition of weak compatibility of splits of~$X$, since the condition on the cardinality is redundant for $k=2$. Condition~\eqref{thm:split-weakly-compatible-3} can never occur if $k=2$, and Condition~\eqref{thm:split-weakly-compatible-2} can only occur in the case $\nu=1$. In this case, however, $i_0,i_3,i_1,i_2\in X$ and the splits $S_1,S_2,S_3$ also fulfil Condition~\eqref{thm:split-weakly-compatible-1} for some $i_0\in X\setminus(S_1(i_1)\cup S_2(i_2) \cup S_3(i_3))$.
\end{enumerate}
\end{proof}
Note that this last proof follows directly from the definition of weak compatibility for splits of sets and splits of polytopes, whereas the proof of \cite[Proposition~6.4]{MR2502496} uses the uniqueness of the split decomposition for metrics \cite[Theorem~2]{MR1153934} and weight functions for polytopes~\cite[Theorem~3.10]{MR2502496}.

We now consider the case $k\geq 3$.

\begin{prop}\label{prop:diss-compatible}
Let $\{A,B\}$, $\{C,D\}$ be two distinct splits of $X$ and $\cT:=\{S_F$ split of $\Hypersimplex kn\,|\,F\in\{A,B,C,D\}\}$ be the set of  corresponding splits of $\Hypersimplex kn$. Then we have:
\begin{enumerate}
\item If $\cT$ is compatible, then $\{A,B\}$ and $\{C,D\}$ are compatible.
\item If $\{A,B\}$ and $\{C,D\}$ are compatible, then there exists at most one non"=compatible pair of splits in $\cT$.
\item If $\{A,B\}$ and $\{C,D\}$ are compatible and $A\cap C=\emptyset$, then $\cT$ is compatible if and only if $k=2$ or $\card{A\cup C}\geq n- k+2$.
\end{enumerate}
\end{prop}
 \begin{proof}
\begin{enumerate}
\item By Lemma~\ref{lem:AB-compatible}, if $\{A,B\}$ and $\{C,D\}$ are not compatible, the only possibility for $S_A$ and $S_C$ or $S_A$ and $S_D$ to be compatible is that $\card{A\cup C}\geq n-k+1$ or $\card{A\cup D}\geq n-k+1$, respectively. However, since $D=X\setminus C$, these two conditions cannot be true at the same time.
\item[(b),(c)] We assume without loss of generality (for (b)) that $A\cut C=\emptyset$. By Lemma~\ref{lem:AB-compatible}, it follows that $S_A$ and $S_B$, $S_B$ and $S_D$, $S_B$ and $S_D$, $S_B$ and $S_C$, and $S_A$ and $S_D$  are compatible, so it only remains to consider the pair $S_A$ and $S_C$. For this pair of splits Lemma~\ref{lem:AB-compatible} implies that it is compatible if and only if $\card{A\cup C} \geq n-k+2$ or $k=2$.
\end{enumerate}
\end{proof}

\begin{cor}\label{cor:compatible-weakly-compatible}
Let $\cS$ be a compatible set of splits of $X$. Then $\cS$ is $k$-weakly compatible for all $k\geq 2$.
\end{cor}
\begin{proof}
This follows directly from Theorem~\ref{thm:split-weakly-compatible}: If either of the properties \eqref{thm:split-weakly-compatible-1}, \eqref{thm:split-weakly-compatible-2}, or \eqref{thm:split-weakly-compatible-3} would hold, then, for example, the pair of splits $\{A_1,X\setminus A_1\}$ and $\{A_2,X\setminus A_2\}$ would not be compatible.
\end{proof}

We conclude by remarking that each of the three conditions in Theorem~\ref{thm:split-weakly-compatible} become weaker as $k$ increases:

\begin{cor}\label{cor:k-l-weakly-compatible}
Let $\cS$ be a set of splits of $X$ and $k\geq3$. If $\cS$ is $k$-weakly compatible, then it is $l$-weakly compatible for all $2\leq l \leq k$. In particular, a $k$-weakly compatible set of splits is weakly compatible.
\end{cor}


\section{$k$"=Dissimilarity Maps from Trees}\label{sec:trees}

Let $T=(V,E,l)$ be a weighted tree consisting of a vertex set $V$, an edge set $E$ and a function $l:E\to \RR_{> 0}$ assigning a weight to each edge. We assume that $T$ does not have any vertices of degree two and that its leaves are labelled by the set~$X$. Such trees are also called \emph{phylogenetic trees}; see Figure~\ref{fig:tree} for an example and Semple and Steel~\cite{MR2060009} for more details. As explained in Figure~\ref{fig:tree}, 
we can define a $k$"=dissimilarity map $D^k_T$ by assigning to each $k$-subset $K\subset X$ the total length of the induced subtree. Each edge $e\in E$ defines a split $S_e=\{A,B\}$ of $X$ by taking as $A$ the set of all leaves on one side of $e$ and as $B$ the set of leaves on the other. It is easily seen that
\begin{align}\label{eq:tree}
D^k_T=\sum_{e\in E}l(e)\delta_{S_e}^k\,.
\end{align}
We now show how this decomposition of $D^k_T$ is related to its 
split decomposition.

\begin{prop}\label{prop:tree-decomposition}
Let $D$ be a $k$"=dissimilarity map on $X$ with $\card X\geq 2k-1$. Then $D=D^k_T$ for some tree $T$ if and only if $\cS_D$ is compatible and $D_0=0$ in the split decomposition of $D$. Moreover, if this holds, then the tree $T$ is unique.
\end{prop}
\begin{proof}
Suppose the split decomposition of $D$ is given by
\[
D=\sum_{\text{$S\in\cS$}} \alpha^D_S \delta^k_S
\]
for some compatible set $\cS$ of splits of $X$. Then Equation~\eqref{eq:tree} shows that for the tree~$T$ whose edges correspond to the splits in $S\in\cS$ with weights $\alpha^D_S$ we have $D^k_T=D$.

Conversely, if $D=D^k_T$ for some weighted tree, Equation~\eqref{eq:tree} is a decomposition of $D^k_T$. By Corollary~\ref{cor:compatible-weakly-compatible}, this decomposition is coherent and the uniqueness part of Theorem~\ref{thm:splitdecomposition} completes the proof. 
\end{proof}

This gives us a new proof of the following Theorem by Pachter and Speyer:

\begin{thm}[\cite{MR2064171}]\label{thm:tree-reconstruct}
Let $T$ be a weighted tree with leaves labelled by $X$ and no vertices of degree two, and $k\geq 2$. If $\card X\geq 2k-1$, then $T$ can be recovered from $D^k_T$.
\end{thm}
\begin{proof}
Compute the split decomposition of $D$. The proof of Proposition~\ref{prop:tree-decomposition} now shows how to construct a tree $T'$ with $D=D^k_{T'}$ and the uniqueness part of this proposition shows that $T=T'$.
\end{proof}

\section{Remarks and Open Questions}

\subsection{Tight-Spans}

\begin{figure}[htb]\centering
  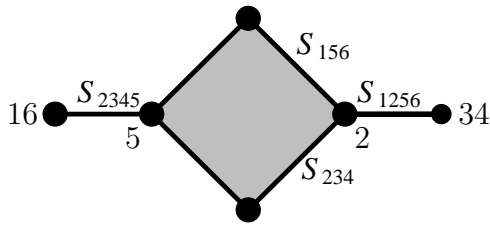
  \caption{The tight-span of the subdivision of $\Hypersimplex 36$ induced by the 3"=dissimilarity map~$D^3_T$ coming from the tree $T$ in Figure~\ref{fig:tree}. Note, that the three non-trivial splits $\{16,2345\}$, $\{34,1256\}$, (corresponding to the splits $S_{2345}$, $S_{1256}$ of $\Hypersimplex 36$, respectively) and $\{156,234\}$ (corresponding to the two splits $S_{156},S_{234}$ of $\Hypersimplex 36$) can be recovered from the tight"=span, as indicated in the figure.}
  \label{fig:tight-span}
\end{figure}

It was shown in \cite[Proposition~2.3]{MR2502496} that the set of inner faces of a regular subdivision $\subdivision wP$ of a polytope~$P$ is anti-isomorphic to a certain realisable polytopal complex, the \emph{tight-span} $\tightspan w P$ of $w$ with respect to $P$. If $P=\Hypersimplex 2n$ and $w_d:=-d$ for a metric $d$ on $X $ then $\tightspan{w_d} {\Hypersimplex 2n}$ is the tight-span $T_d$ of the metric space $(X,d)$; see Isbell~\cite{MR82949} and Dress~\cite{MR753872}. In particular, if $d$ is a tree metric, then $T_d$ is isomorphic to that tree. For a $k$"=dissimilarity map~$D$ one can similarly consider the tight-span  $\tightspan{w_D} {\Hypersimplex kn}$. However, Proposition~\ref{prop:diss-compatible} shows that $\tightspan{w_D} {\Hypersimplex kn}$ is not necessarily a tree for $k\geq 3$. As an example, we depict in Figure~\ref{fig:tight-span} the tight-span $\tightspan{w_{D^3_T}}{\Hypersimplex 36}$ where $T$ is the tree from Figure~\ref{fig:tree}. Even though it is not a tree, note that the non"=trivial splits corresponding to the edges of $T$ can be easily recovered from $\tightspan{w_{D^3_T}}{\Hypersimplex 36}$. It would be interesting to understand better the relationship between the structure of $\tightspan{w_D}{\Hypersimplex kn}$ and the split decomposition of~$D$ in case $D$ has no split-prime component.

\subsection{Matroid Subdivisions, Tropical Geometry, and Valuated Matroids}
A subdivision $\subdiv$  of $\Hypersimplex kn$ is called a \emph{matroid subdivision} if all $1$-dimensional cells $E\in\subdiv$ are edges of $\Hypersimplex kn$, or, equivalently, if all elements of $\subdiv$ are matroid polytopes. The space of all weight functions $w$ inducing matroid subdivisions is called the \emph{Dressian}. The elements of the Dressian correspond to (uniform) valuated matroids (see \cite[Remark~2.4]{MR2515769}) and to tropical Pl\"ucker vectors (see Speyer~\cite[Proposition~2.2]{MR2448909}). The corresponding weight function $w$ then defines a so called \emph{matroid subdivision} of $\Hypersimplex kn $. The \emph{tropical Grassmannian} (see \cite{MR2071813}) is a subset of the Dressian. It was shown by Iriarte \cite{Iriarte} with methods developed by Bocci and Cools \cite{MR2531242}, and   Cools \cite{MR2523769} that for a weighted tree~$T$, the weight function $w_{D^k_T}$ is a point in the tropical Grassmannian and hence in the Dressian. Corollary~\ref{cor:compatible-weakly-compatible} now implies that $w_{D^k_T}$ is indeed in the interior of the cone of the Dressian spanned by the split weights~$w^k_{S_e}$ for all splits $S_e$ corresponding to edges $e$ of~$T$. In the language of matroid subdivisions this implies that starting from a compatible set~$\cS$ of splits of~$X$ the set $\smallSetOf{S_A\text{ split of }\Hypersimplex kn}{A\in S, S\in\cS}$ of splits of $\Hypersimplex kn$ induces a matroid subdivision. Establishing that other sets of splits satisfying the requirements of Theorem~\ref{thm:hs-weakly-compatible} also have this property could lead to a further understanding of the Dressian.

\subsection{Computation of the Split Decomposition and Tree Testing}
In \cite{MR2064171}, Speyer and Pachter raise the question how to test whether a given $k$"=dissimilarity map $D$ on $X$ comes from a tree. Our results suggest the following simple algorithm: Compute the split indices $\alpha^D_S$ for all splits of~$X$, test whether $D_0=0$ in the split decomposition~\eqref{eq:split-decomposition}, and whether the split system $\cS_D$ is compatible. Equation~(2) in~\cite{MR2502496} gives an explicit formula for the indices $\alpha^{w_D}_{w_{S_A}}$ and hence for the split indices $\alpha^D_S$, however this involves the computation of the tight-span $\tightspan{w_D} {\Hypersimplex kn}$ whose number of vertices can be in general exponential in $n$. It would be interesting to derive a simpler formula for the split indices similar to the one existing in the case $k=2$ given by Bandelt and Dress~\cite[Page~50]{MR1153934}. This might yield a polynomial algorithm to test whether a given $k$"=dissimilarity map $D$ on $X$ comes from a tree.

\bibliographystyle{amsplain-mod}
\bibliography{tight_spans}

\end{document}